\documentclass[10pt]{amsart}
\usepackage{mathrsfs}
\usepackage{amstext}
\usepackage{amscd}
\usepackage{latexsym}
\usepackage{amsfonts}
\usepackage{amssymb}
\usepackage{amsmath}
\usepackage{amsthm}
\usepackage{color}
\usepackage{multicol}
\usepackage{pgfplots}
\usepackage{multirow}
\usepackage{longtable}

\usepackage[cp1251]{inputenc}
		\theoremstyle{plain}
		\newtheorem{theorem}{Theorem}[section]
		\newtheorem{lemma}[theorem]{Lemma}
		\newtheorem{prop}[theorem]{Proposition}
		
		\newtheorem{conjecture}[theorem]{Conjecture}		
		
		\theoremstyle{definition}
		\newtheorem{definition}[theorem]{Definition}
		
		\newtheorem{remark}[theorem]{Remark}

		\newtheorem{problem}[theorem]{Problem}
		
		\pgfplotsset{compat=1.16}

\begin{document}

\title[Integral Laplacian graphs with a unique double Laplacian eigenvalue]{Integral Laplacian graphs with a unique double Laplacian eigenvalue, I}			
\author[A.~Hameed]{Abdul Hameed}
			
\address{School of Mathematical Sciences, Shanghai Jiao Tong University, Shanghai, P.R.China}
\email{abdul-hameed\_211@sjtu.edu.cn}
			
\author[M.~Tyaglov]{Mikhail Tyaglov}
\address{School of Mathematical Sciences and MOE-LSC, Shanghai Jiao Tong University, Shanghai, P.R.China}
\email{tyaglov@gmail.com, tyaglov@sjtu.edu.cn}


\begin{abstract}
The set $S_{i,n}=\{0,1,2,\ldots,n-1,n\}\setminus\{i\}$, $1\leqslant i\leqslant n$ is called Laplacian realizable if there exists an undirected simple graph
whose Laplacian spectrum is $S_{i,n}$. The existence of such graphs was established by S.\,Fallat et al. in 2005.
In this paper, we investigate graphs whose Laplacian spectra have the form
\begin{equation*}
S_{\{i,j\}_{n}^{m}}=\{0,1,2,\ldots,m-1,m,m,m+1,\ldots,n-1,n\}\setminus\{i,j\},\qquad 0<i<j\leqslant n,
\end{equation*}
and completely describe those ones with $m=n-1$ and $m=n$. We also show close relations between graphs
realizing~$S_{i,n}$ and $S_{\{i,j\}_{n}^{m}}$, and discuss the so-called $S_{n,n}$-conjecture and the 
correspondent conjecture~for~$S_{\{i,n\}_{n}^{m}}$.
\end{abstract}

\maketitle
			
\setcounter{equation}{0}
			


\setcounter{equation}{0}
\section{Introduction}\label{Section:Introduction}

Let $G=(V(G),E(G))$ be an undirected simple graph with vertex set $V(G)={\{v_{1},v_{2},\ldots,v_{n}}\}$ and the edge set $E(G)=\{e_{1},\,e_{2},\ldots,\,e_{r}\}$, where $\left|V(G)\right|=n$ and $\left|E(G)\right|=r$ are respectively called the order and the size of~$G$.  We denote the degree of vertex $v_{i}$, $i=1,2,\ldots,n$ by $d_{i}=d(v_{i})$. The \textit{Laplacian matrix} of~$G$ is a matrix defined  as follows: $L{(G)}=D{(G)}-A{(G)}$, where $D{(G)}$ is the vertex degree diagonal matrix, i.e., $D({G})=\mathrm{diag}(d_{1},\,d_{2},\ldots,\,d_{n-1},\,d_{n})$ and $A{(G)}$ is the $(0,1)$ adjacency matrix of~$G$. Thus, the entries of the Laplacian matrix have the following form

\begin{equation*}
l_{ij} =
\begin{cases}
&\ \,d_i,\quad\text{if}\quad \ \ i=j,\\
&-1,\ \ \ \text{if}\quad \ \  i\neq j\ \text{and}\  v_{i}\sim v_{j}, \\
&\ \ \,0,\ \ \ \text{otherwise}.
\end{cases}
\end{equation*}

\vspace{4mm}

The study of graphs whose adjacency matrix has integer eigenvalues was probably introduced by F.\,Harary and A.J.\,Schwenk~\cite{HararySchwenk_1974}. The same kind of problem has been addressed to the eigenvalues of the Laplacian matrix. A graph $G$ whose Laplacian matrix has integer eigenvalues is called \textit{Laplacian integral}. Some classes of Laplacian integral graphs have been identified. For example, R.\,Merris~\cite{Merris.1_1994} showed that degree maximal graphs are  Laplacian integral. J.J.\,Kirkland~\cite{Kirkland_2008} characterized Laplacian integral graphs of maximum degree~$3$.   The most well-known examples for such graphs are the complete graph $K_{n}$, the star graph $S_{n}$, and the Petersen graph having eigenvalues of its Laplacian matrix $\{0,{n}^{n-1}\}$, $\{0,{1}^{n-2},n\}$ and $\{0,2^{5},5^{4}\}$, respectively (the superscript indicates the multiplicity of the eigenvalue). Every \textit{threshold graph} (i.e., a graph that can be constructed from isolated vertices by the operations of unions and complements) is Laplacian integral, as well. For more detail on Laplacian integral graphs, we refer the reader to the 
works~\cite{Balinska_2002,Lima_et_al_2007,FallatKirkland_et_al_2005,GroneMerris_2008,HammerKelmans_1996,Kirkland_2005,Kirkland_2007,Kirkland_2008,Kirkland_et_al_2010,Merris.1_1994,Merris_1997,Merris_1997.1}.

One of the way to describe and study all the integral Laplacian graphs (if any) is to study graphs having certain kind of Laplacian spectra.
Thus, starting with graphs with simple integer Laplacian eigenvalues, one then can consider integral Laplacian graphs with all simple Laplacian eigenvalues
excluding one which is double. On the next step, it is allowed for graphs to have two multiple Laplacian eigenvalues etc.

The first step in this scheme was undertaken by S.\,Fallat et al. in their work~\cite{FallatKirkland_et_al_2005} on the integral Laplacian graphs with simple Laplacian eigenvalues.
\begin{definition}
A set $S$ consisting of elements of $\{0,1,\ldots,n\}$ (maybe repeated) such that $|S|=n$ and $0\in S$ is called Laplacian realizable if there exists a connected
simple graph on $n$ vertices whose Laplacian spectrum coincides with $S$. In that situation, the graph $G$ is said to
realize~$S$.
\end{definition}
Thus, the authors in~\cite{FallatKirkland_et_al_2005} considered the sets of the form
\begin{equation}\label{set_S_i,n}
S_{i,n}=\{0,1,2,\ldots,n-1,n\}\setminus\{i\},
\end{equation}
and found out when such sets are Laplacian realizable and completely described the graphs realizing $S_{i,n}$. The only problem
that remained open in this case is whether the set $S_{n,n}$ is Laplacian realizable or not. This problem is now known
as the $S_{n,n}$-conjecture and states that $S_{n,n}$ is \textit{not} Laplacian realizable for every $n\geqslant2$. This conjecture
was proved for $n\leqslant11$, for prime $n$, and for $n\equiv 2,3\mod4$ in~\cite{FallatKirkland_et_al_2005}. Later, Goldberger
and Neumann~\cite{GoldbergerNeumann_2013} showed that the conjecture is true for $n\geqslant6,649,688,933$. The authors of the present
work established~\cite{AhZuMt_2021} that if a graph is the Cartesian product of two other graphs, then it does not realize $S_{n,n}$.

\vspace{2mm}

Following the aforementioned scheme, in the present paper, we study Laplacian integral graphs having only one multiple (non-zero) Laplacian eigenvalue whose
multiplicity equals exactly $2$. Namely, we describe the graphs realizing the sets of the form
\begin{equation}\label{set_S_i.j.n}
S_{\{i,j\}_{n}^{m}}=\{0,1,2,\ldots,m-1,m,m,m+1,\ldots,n-1,n\}\setminus\{i,j\},
\end{equation}
where $0<i<j\leqslant n$. So $S_{\{i,j\}_{n}^{m}}$ does not contain the numbers $i$ and $j$, while some number $m$ (and only this number)
is doubled. Obviously, if $m=0$, then the graph is disconnected, so we exclude this case from the further considerations.

\begin{figure}[h]
\center{\includegraphics[width=9cm,height=3cm]{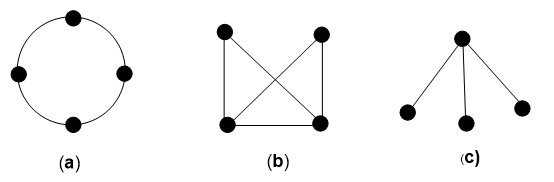}}\caption{Graphs of order $4$ realizing $S_{\{1,3\}_{4}^{2}}$, $S_
{\{1,3\}_{4}^{4}}$ and $S_{\{2,3\}_{4}^{1}}$ respectively.}\label{order4.pic}
\centering
\end{figure}

The graphs realizing $S_{\{i,j\}_{n}^{m}}$ exist even for small $n$. In~\cite[p.~286--289]{CvetkovicRowlinson_2010} the authors
found the Laplacian spectra of all graphs of order up to $5$, so it easy to observe that for $n=3$, only the set $S_{\{1,2\}_{3}^{3}}$
is Laplacian realizable, and it is realized by the complete graph $K_3$. For $n=4$, there are exactly three sets $S_{\{1,3\}_{4}^{2}}$,
$S_{\{1,3\}_{4}^{4}}$ and~$S_{\{2,3\}_{4}^{1}}$ that are Laplacian realizable. They are depicted on Figure~\ref{order4.pic}.
For $n=5$, the only Laplacian realizable sets of the form $S_{\{i,j\}_{n}^{m}}$ are $S_{\{2,3\}_{5}^{4}}$, $S_{\{1,3\}_{5}^{5}}$,
$S_{\{2,4\}_{5}^{3}}$, $S_{\{1,4\}_{5}^{2}}$ and $S_{\{2,4\}_{5}^{1}}$.

Note that the graphs $(b)$ and $(c)$ in Figure~\ref{order4.pic} are threshold, see, e.g.,~\cite{Merris.1_1994,Merris_1997.1}. The explicit constructions of graphs
realizing~$S_{\{i,j\}_{n}^{m}}$ for $n=4,5,6$ are given in Tables~\ref{Table.2}, \ref{Table.3}, \ref{Table.4} of Appendix~\ref{Tables.constructions}.
To make these tables we used the complete lists of all connected graphs up to order $6$, see~\cite[p.~286--289]{CvetkovicRowlinson_2010} and~\cite{CvetkovicPetric_1984}

The present work is the first part of our research on graphs realizing the sets $S_{\{i,j\}_{n}^{m}}$, and is devoted to the cases $m=n-1$ and $m=n$.
We list all Laplacian realizable sets $S_{\{i,j\}^n_n}$, Theorem~\ref{Case_m=n}, and describe the structure of the graphs realizing
$S_{\{i,j\}^n_n}$, Theorem~\ref{Constru.Double.Eig}. For the case $m=n-1$, we show that only~$S_{\{1,j\}_n^{n-1}}$
and~$S_{\{2,j\}_n^{n-1}}$ can be Laplacian realizable for certain $j$, Theorems~\ref{condition.n-1}, and list all such $j$ for a given~$n$,
Theorems~\ref{double.m=n-1} and~\ref{Theorem.S.1.j.n-1.n_constr}. Theorems~\ref{Thm.case.2.m=n-1} and~\ref{Theorem.S.1.j.n-1.n} describe the structure of graphs realizing
the sets~$S_{\{i,j\}_n^{n-1}}$ for~$i=1,2$.

\begin{figure}[h]
\begin{center}
		
\includegraphics[width=7.5cm,height=2.5cm]{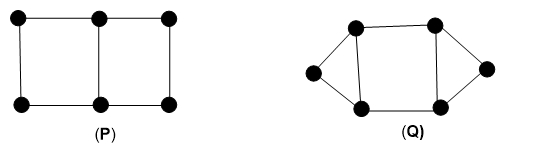}
\caption{Graph $P$ and its complement $Q$ realizing $S_{\{4,6\}_6^3}$ and $S_{\{2,6\}_6^3}$, respectively.}\label{ladder1.pic}
\centering
\end{center}
\end{figure}

\begin{figure}[h]
\begin{center}
\includegraphics[width=8cm,height=3.5cm]{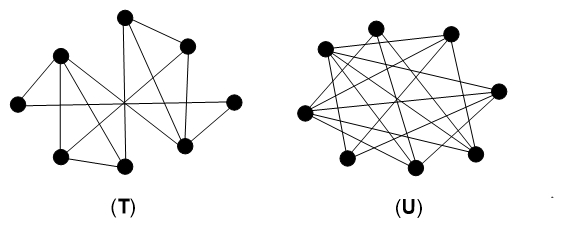}
\caption{Graph $T$ and its complement $U$ realizing $S_{\{7,8\}_8^3}$ and  $S_{\{1,8\}_8^5}$, respectively }\label{ladder2.pic}
\centering
\end{center}
\end{figure}

Unlike the sets $S_{n,n}$ which are conjectured to be not Laplacian realizable, the sets $S_{\{i,n\}_{n}^{m}}$ (that is, without $n$ like $S_{n,n}$)
can be Laplacian realizable, at least for small $n$. In Section~\ref{section.conj}, we show that the ladder graphs on~$6$ and a graph on $8$ vertices and their complement realize such sets (see Figure~\ref{ladder1.pic} and \ref{ladder2.pic}).

\noindent However, there are no other Laplacian realizable sets $S_{\{i,n\}_n^m}$ for $n\leqslant7$, and we conjecture that there are no
such sets for $n\geqslant6$ except the 4 graphs depicted in Figures~\ref{ladder1.pic}--\ref{ladder2.pic}, see Conjecture~\ref{Conjecture.S_nnm}.
In attempts to prove this conjecture, we found out that the sets~$S_{\{i,n\}_n^m}$
are not Laplacian realizable for prime integer $n$, Theorem~\ref{conj.facts}, and that any graph realizing~$S_{\{i,n\}_n^m}$ (if any)
cannot be the Cartesian product of two graphs for $n\geqslant9$, Theorem~\ref{Cartesian.Prod.j=n}. To establish some results of the present work, we used new
Laplacian spectral properties of the join of graphs, Theorems~\ref{corol.Join.n.n} and~\ref{eig.1.join}.

The paper is organized as follows. In Section~\ref{section:preliminaries}, we list basic definitions and important results that we use in our
work. Also we prove some auxiliary statements and theorems. 
In Sections~\ref{section.results.1} and~\ref{section.results.2}, we list
all the Laplacian realizable sets of the form $S_{\{i,j\}_{n}^{n}}$ and $S_{\{i,j\}_{n}^{n-1}}$, respectively, and develop
an algorithm for constructing the graphs realizing those sets. Section~\ref{section.conj} is devoted to the Conjecture~\ref{Conjecture.S_nnm}. Some concluding
remarks are presented in Section~\ref{section:conclusion}. Finally, in Appendix~\ref{Tables.constructions},
we list all the Laplacian realizable sets $S_{\{i,j\}_{n}^{m}}$ for $n=4,5,6$ and a few ones for $n=7$. The correspondent
graphs realizing those sets are presented.

\setcounter{equation}{0}

\section{Preliminaries}\label{section:preliminaries}

Let $G$ be an undirected simple graph. The \textit{complement} $\overline G$ of the graph $G$ is a graph $\overline{G}$ on the same vertices such that two distinct vertices of $\overline{G}$ are adjacent if and only if they are not adjacent in $G$.

Given two graphs $G_{1}=(V_{1}, E_{1})$ and~$G_{2}=(V_{2}, E_{2})$, the \textit{union} of graphs $G_{1}$ and $G_{2}$, denoted as $G_{1}\cup G_{2}$, is the graph $G=(V,E)$ for which $V=V_{1}\cup V_{2}$ and $E=E_{1}\cup E_{2}$. The union of $k$ copies of the same graph $G$ is denoted as $kG$. The \textit{join} $G_{1}\vee G_{2}$ of the graphs~$G_{1}$ and~$G_{2}$ is the graph obtained from $G_{1}\cup G_{2}$ by joining every vertex of $G_{1}$ with every vertex~of~$G_{2}$, that is, $G_{1}\vee G_{2}=
\overline{(\overline{G_{1}}\cup \overline{G_{2}})}$.

The \textit{Cartesian product} of the graphs $G_1$ and $G_2$ is the graph $G_1\times G_2$ whose vertex set is the Cartesian product $V(G_1)\times V(G_2)$, and for $v_{1},v_{2}\in V(F)$ and  $u_{1},u_{2}\in V(H)$, the vertices $(v_{1},u_{1})$ and $(v_{2},u_{2})$ are adjacent in $G_1\times G_2$ if and only if either
\begin{itemize}
\item[$\bullet$] $v_{1}=v_{2}$ and $\{u_{1},u_{2}\}\in E(G_1)$;
\item[$\bullet$]${\{v_{1},v_{2}\}}\in E(G_2)$ and $u_{1}$=$u_{2}$.
\end{itemize}	
One of the examples of the Cartesian product is the ladder graph which is the Cartesian product of two path graphs, one of which has only one edge: $L_{n,1}=P_n\times P_2$.

We denote the eigenvalues of the Laplacian matrix $L(G)$ arranged in increasing order as	
$$0=\mu_{1}\leqslant\mu_{2}\leqslant\ldots\leqslant\mu_{n-1}\leqslant\mu_{n}.$$ It is well
known that $L{(G)}$ has a zero eigenvalue corresponding to the eigenvector with equal entries,
while other eigenvalues are non-negative~\cite{Merris_1998}. M.\,Fiedler~\cite[p.~298]{Fiedler_1973}
showed that $\mu_{2}>0$ if and only if the graph~$G$ is connected. This eigenvalue is usually known as
the \textit{algebraic connectivity} of~$G$ denoted by $a_{G}$. The \textit{Laplacian spectrum}
$\sigma_L(G)$ of a graph $G$ is the spectrum of the Laplacian matrix of~$G$. The \textit{Laplacian
spectral radius} denoted as $\rho(G)$ of~$G$ is the absolute value of the largest eigenvalue of
the Laplacian matrix~$L{(G)}$, i.e., $\rho(G)=\max\limits_{1\leqslant i\leqslant n}\left|\mu_{i}\right|$.
The Laplacian spectrum of a graph is bounded from above by the order of the graph, see,~\cite{Merris_1994}.
\begin{prop}\label{Prop.max.eig}
Let $G$ be a simple graph on $n$ vertices. Then $\rho(G)\leqslant n$.
\end{prop}

The Laplacian spectrum of the graph operations mentioned above is related to the Laplacian spectra of the initial graphs.
Thus, the Laplacian spectrum of $\bigcup\limits_{j=1}^kG_{j}$ is the union of the Laplacian spectra of the graphs
$G_{1}$, \ldots, $G_{k}$, see, e.g.,~\cite{CvetkovicRowlinson_2010}. The Laplacian spectra of the complement of a graph
and the join of two graphs are given by the following theorems, see e.g.~\cite{CvetkovicRowlinson_2010,Merris_1994}.
\begin{theorem}\label{Complement.spect}
Let $G$ be a graph with $n$ vertices. If $\sigma_{L}(G)=\{0,\mu_2,\ldots,\mu_{n-1},\mu_n\}$, then
\begin{equation*}
\sigma_{L}(\overline{G})=\{0,n-\mu_n,n-\mu_{n-1},\ldots,n-\mu_2\}.
\end{equation*}
\end{theorem}

The spectrum of the join of two graphs was obtain by Kel'mans~\cite{Kelmans_1965} in terms of the characteristic polynomials of Laplacian matrices. The following form of Kelmans' theorem can be found, e.g., in~\cite{Merris_1998}.
\begin{theorem}\label{Join.Thm}
Let $G$ and $H$ be two graphs of order $n$ and $m$, respectively, whose eigenvalues are the following
\begin{equation*}
0=\mu_{1}\leqslant\mu_{2}\leqslant\mu_{3}\ldots\leqslant \mu_{n-1}\leqslant\mu_{n} \qquad\text{and}\qquad 0=\lambda_{1}\leqslant\lambda_{2}\leqslant\lambda_{3}\ldots\leqslant\lambda_{m-1}\leqslant\lambda_{m}.
\end{equation*}
Then the Laplacian spectrum of the join $G\vee H$ has the form
\begin{equation}\label{join.spectrum.1}
\{0,m+\mu_{2},m+\mu_{3},\ldots,m+\mu_{n-1},m+\mu_{n},n+\lambda_{2},n+\lambda_{3},\ldots,n+\lambda_{m-1},n+\lambda_{m}, n+m\}.
\end{equation}
\end{theorem}

\noindent Note that the eigenvalues here are not in increasing order.

\vspace{2mm}

As one can see from~\eqref{join.spectrum.1}, the order of the join of two graphs is a Laplacian eigenvalue of the join. It turns out that this fact is a necessary and sufficient condition for a graph to be a join, see, e.g.,~\cite{Molitierno_2016}.
\begin{theorem}\label{Thm.Join.n}
Let $G$ be a connected graph of order $n$. Then $n$ is a Laplacian eigenvalue of~$G$ if and only if $G$ is the join of two graphs.
\end{theorem}

However, if the order of a graph is a Laplacian eigenvalue of multiplicity $2$, then we can get more information about the graph.

\begin{theorem}\label{corol.Join.n.n}
Let $G$ be a connected graph of order~$n$, and let $n$ be the Laplacian eigenvalue of~$G$ of multiplicity~$2$. Then $G=F\vee H$ where $F$ is a join of two graphs, while $H$ is not a join.
Moreover, the number~$1$ does not belong to the Laplacian spectrum of~$G$ in this case.
\end{theorem}
Thus, the algebraic connectivity of any graph of order $n$ with doubled Laplacian eigenvalue $n$ is greater than~$1$.
\begin{proof}
Indeed, since $n$ is a Laplacian eigenvalue of~$G$, the graph $G$ is a join of two graphs according to  Theorem~\ref{Thm.Join.n}, that is, $G=F\vee H$.
Suppose that $\left|F \right|=p$ and $\left|H \right|=n-p$ for some integer number $p$, $1\leqslant p\leqslant n-1$.
	
Let us denote the Laplacian spectra of the graphs $F$ and $H$ as follows
\begin{equation}\label{EV.G.mu}
0=\mu_{1}\leqslant\mu_{2}\ldots\leqslant \mu_{p-1}\leqslant\mu_{p}\quad\text{and}\quad 0=\lambda_{1}\leqslant\lambda_{2}\ldots\leqslant\lambda_{n-p-1}\leqslant\lambda_{n-p}.
\end{equation}
Then by Theorem~\ref{Join.Thm}, the Laplacian spectrum of~$G$ has the form
\begin{equation}\label{spectrum.join}			
\sigma_{L}(G)=\{0,(n-p)+\mu_{2},\ldots,(n-p)+\mu_{p},p+\lambda_{2},\ldots,p+\lambda_{n-p-1},p+\lambda_{n-p},n\}.
\end{equation}
Here the eigenvalues are not in the increasing order.
	
Since the graph $G$ has exactly two eigenvalues $n$ by assumption, one of the numbers $n-p+\mu_{p}$ or $p+\lambda_{n-p}$ equals $n$ while the other one is less than $n$. Without loss of generality, we can suppose that $n-p+\mu_{p}=n$ and $p+\lambda_{n-p}<n$ that implies $\mu_{p}=p$ and $\lambda_{n-p}<n-p$.
	
Now from~\eqref{spectrum.join} we obtain that $\mu_{p}=p$ is a simple Laplacian eigenvalue of~$G$. Indeed, if $\mu_{p}$ is not simple, then $n$ has multiplicity at least $3$ that contradicts the assumption of the theorem.
Consequently, the graph $F$ is a join with simple maximal Laplacian eigenvalue, whereas the graph $H$ is not a join by Theorem~\ref{Thm.Join.n}, since $\lambda_{n-p}<n-p$.
	
According to~\eqref{EV.G.mu}--\eqref{spectrum.join}, either (or both) numbers $(n-p)+\mu_{2}$ and $p+\lambda_{2}$ is the minimal positive Laplacian eigenvalue of the graph $G$. However, since the graph $F$ of order $p$ is a join, it is clear that $p\geqslant2$, at the same time $n-p\geqslant1$. So both numbers $(n-p)+\mu_{2}$ and $p+\lambda_{2}$ are greater than $1$, as required.
\end{proof}

The following proposition provides a necessary and sufficient condition on a graph $G$ to have a Laplacian eigenvalue $1$ provided $G$ is a join.
\begin{theorem}\label{eig.1.join}
Let a graph $G$ be a join. The number $1$ is a Laplacian eigenvalue of~$G$ if and only if $G=F\vee K_{1}$ where~$F$ is a disconnected graph of order at least $2$.
\end{theorem}
\begin{proof}
Let $G$ be a join that has a Laplacian eigenvalue $1$ and let $G=F\vee H$ where the graph $F$ is of order~$p$ while the graph $H$ is of order $n-p$. According to~\eqref{EV.G.mu}--\eqref{spectrum.join}, one of the numbers $(n-p)+\mu_{2}$ or $p+\lambda_{2}$ equals~$1$.
	
If both of the numbers $(n-p)+\mu_{2}$ and $p+\lambda_{2}$ equal $1$, then $\mu_2=\lambda_2=0$ and $(n-p)=p=1$. So in this case $G$ is a join of two isolated vertices, i.e., $G=K_{1}\vee K_{1}$, and the Laplacian spectrum of~$G$ is equal to $\{0,2\}$. But this contradicts to the assumption that $G$ has a Laplacian eigenvalue $1$.
	
Thus, we have that one of the numbers $(n-p)+\mu_{2}$ and $p+\lambda_{2}$ is greater than $1$. Without loss of generality, we can suppose that $(n-p)+\mu_{2}=1$ and $p+\lambda_{2}>1$ (as $p$ is a positive integer number). This implies that $(n-p)=1$ and $\mu_{2}=0$. So the graph $F$ has at least two Laplacian eigenvalues equal zero. Consequently, $p\geqslant2$ and $F$ is disconnected~\cite[p.~298]{Fiedler_1973}. As $(n-p)=1$, the graph $H$ consists of a single vertex~$K_{1}$.

	
Conversely, let $G=F\vee K_{1}$ where $F$ is a disconnected graph of order at least $2$. Since $F$ is disconnected, the multiplicity of the Laplacian eigenvalue $0$ is at least $2$, so by Theorem~\ref{Join.Thm}, the Laplacian spectrum of~$G$ contains $1$ is an eigenvalue of~$G$.
\end{proof}

The Laplacian spectrum of the Cartesian product of two graphs is given in the following theorem, see,~e.g.,~\cite{Merris_1994}.
\begin{theorem}\label{Cartesian.eigen}
Let $F$ and $H$ be graphs having spectrum,
\begin{equation*}
	\sigma_{L}(F)=(\lambda_{1},\lambda_{2},\ldots,\lambda_{n-1},\lambda_{n}),\;\; \sigma_{L}(H)=(\mu_{1},\mu_{2},\ldots,\mu_{p-1},\mu_{p}),
\end{equation*}
then the Laplacian spectrum of the Cartesian product of $F$ and $H$ is
\begin{equation}\label{Cartesian.Prod.Spectrum}
\sigma_{L}(F\times H)=\{\lambda_{i}+\mu_{k}\},\quad 1\leqslant i\leqslant n, \quad 1\leqslant k\leqslant p.
\end{equation}
\end{theorem}
From~\ref{Cartesian.Prod.Spectrum} it is easy to see that the Laplacian spectrum of the Cartesian product of two Laplacian integral graphs is Laplacian integral.


As we mentioned in Introduction, the graphs whose Laplacian spectra are of the form~$S_{i,n}$ defined in~\eqref{set_S_i,n}, were
introduced and completely investigated in~\cite{FallatKirkland_et_al_2005}. In the sequel, we use the following their results~\cite[Theorems~2.3 and 2.6]{FallatKirkland_et_al_2005}.

\begin{prop}\label{All.Sin}
%
\begin{itemize}
\item[(i)] If  $n\equiv 0\mod 4$, then for each $ i=1,2,3,\dots,\frac{n-2}{2}$, $S_{2i,n}$ is Laplacian realizable;
\item[(ii)] If $n\equiv 1\mod 4$, then for each $ i=1,2,3,\dots,\frac{n-1}{2}$, $S_{2i-1,n}$ is Laplacian realizable;
\item[(iii)] If $n\equiv 2\mod 4$, then for each $ i=1,2,3,\dots,\frac{n}{2}$, $S_{2i-1,n}$ is Laplacian realizable;
\item[(iv)] If $n\equiv 3\mod 4$, then for $ i=1,2,\dots,\frac{n-1}{2}$, $S_{2i,n}$ is Laplacian realizable.
\end{itemize}
\end{prop}
%
%
%
%
%
%
%
%
%
%
%
%
\begin{prop}\label{Propos:Constr.S_i.n.222}
Suppose that $n\geqslant 6$ and that $G$ is a graph on $n$ vertices. Then $G$ realizes~$S_{1,n}$ if and only if $G$ is formed in
one of the following two ways:
\begin{itemize}
\item[(i)] $G=(K_{1}\cup K_{1})\vee (K_{1}\cup G_{1})$, where  $G_{1}$ is a graph on $n-3$ vertices that realizes $S_{n-4,n-3};$
\item[(ii)] $G=K_{1}\vee H$, where  $H$ is a graph on $n-1$ vertices that realizes $S_{n-1,n-1}$.
\end{itemize}
\end{prop}

\begin{prop}\label{Propos:Constr.S_i.n.223}
Suppose that $n\geqslant 6$ and that $G$ is a graph on $n$ vertices. Then $G$ realizes~$S_{i,n}$ with $2\leqslant i\leqslant n-2$
if and only if $G=K_{1}\vee (K_{1}\cup H)$, where $H$ is a graph on $n-2$ vertices that realizes $S_{i-1,n-2}$.
\end{prop}
\begin{prop}\label{Propos:Constr.S_i.n.224}
Suppose that $n\geqslant 6$ and that $G$ is a graph on $n$ vertices.
$G$ realizes $S_{n-1,n}$ if and only if~$G$ is formed in one of the following two ways:
\begin{itemize}					
\item[(i)]  $G=K_{1}\vee (K_{2}\cup G_{1})$, where $G_{1}$ is a graph on $n-3$ vertices that realizes $S_{2,n-3}$
\item[(ii)] $G=K_{1}\vee (K_{1}\cup H)$, where $H$ is a graph on $n-2$ vertices that realizes $S_{n-2,n-2}$.
\end{itemize}
\end{prop}
Proposition~\ref{All.Sin} completely resolves the existence of graphs realizing
the spectrum~$S_{i,n}$, where $1\leqslant i<n$. Furthermore, it is easily deduced from
Propositions~\ref{Propos:Constr.S_i.n.222}--\ref{Propos:Constr.S_i.n.224} that if the set $S_{n,n}$
were not realizable for any $n$, then there is a unique graph, which realizes $S_{i,n}$, for $1\leqslant i<n$.
As we mentioned in Introduction, it is already proved that for $n\leqslant11$ and $n\geqslant6,649,688,933$, the sets
$S_{n,n}$ are not Laplacian realizable, see~\cite{FallatKirkland_et_al_2005,GoldbergerNeumann_2013}.
Below in Sections~\ref{section.results.1}--\ref{section.results.2} we prove some analogues of these results for
the sets $S_{\{i,j\}_{n}^{n}}$ and $S_{\{i,j\}_{n}^{n-1}}$.

\vspace{2mm}

Now we are in a position to establish a few general facts on the realizability of sets $S_{\{i,j\}_{n}^{m}}$.
\begin{theorem}\label{Condition.m}
Suppose that $n\geqslant 3$ and $G$ is a graph of order $n$ realizing $S_{\{i,j\}_{n}^{m}}$ for $i<j$. Then
\begin{itemize}
\item[(i)] for $n\equiv 0\;or\; 3 \mod 4$, the numbers $(i+j)$ and $m$ are of the same parity;
\item[(ii)] for $n\equiv 1\;or\;2 \mod 4$, the numbers $(i+j)$ and $m$ are of opposite parity.
\end{itemize}
\end{theorem}
\begin{proof}
If a graph $G$ realizes $S_{\{i,j\}_{n}^{m}}$, then the sum of elements in the set $S_{\{i,j\}_{n}^{m}}$ equals the sum of degrees of vertices in the graph $G$ (the trace of $L(G)$). So
\begin{equation}\label{main.equation}
\mathrm{Tr}\,(L(G))=2\left|E(G)\right|=\dfrac{n(n+1)}{2}+m-(i+j).
\end{equation}
Here $m$ is the double eigenvalue of $L(G)$.

Now if $n\equiv 0$ or $3 \mod 4$, then $\frac{n(n+1)}{2}$ is an even number. Since $\mathrm{Tr}\,(L(G))$ is even as well,
one can see that~$m-(i+j)$ must be even by~\eqref{main.equation}. Therefore, $m$ and $i+j$ are of the same parity. Analogously,
if $n\equiv 1$ or $2 \mod 4$, then $\frac{n(n+1)}{2}$ is odd, so $(i+j)$ and $m$ are of the opposite parity, as required.
\end{proof}

The following lemma is of frequent use in the sequel.
\begin{lemma}\label{join.union}
	Let $n\geqslant 3$, if $S_{\{i,j\}_{n}^{m}}$ is Laplacian realizable, then so is $S_{\{i+1,j+1\}_{n+2}^{m+1}}$.
\end{lemma}
\begin{proof}
By Theorem~\ref{Join.Thm}, one easily obtains that if a graph $G$ realizes the set $S_{\{i,j\}_{n}^{m}}$, then the
Laplacian spectrum of the graph $K_{1}\vee(K_{1}\cup G)$ is exactly $S_{\{i+1,j+1\}_{n+2}^{m+1}}$.
\end{proof}

\begin{remark}
Note that the converse statement of Lemma~\ref{join.union} does not hold, in general. For instance, from Table~\ref{Table.3}
in Appendix~\ref{Tables.constructions} it follows that the set~$S_{\{2,4 \}_5^3}$ is Laplacian realizable. However,
$S_{\{1,3 \}_3^2}$ is not Laplacian realizable, since the only realizable set $S_{\{i,j \}_3^m}$ is $S_{\{1,2\}_{3}^{3}}$ as
we mentioned in Introduction. Analogously, $S_{\{4,6 \}_6^3}$ is Laplacian realizable, but $S_{\{3,5\}_4^2}$ is not.
\end{remark}

We conclude this section with the following theorem showing that for large $n$ the graphs realizing $S_{\{i,j\}_{n}^{m}}$
cannot be the Cartesian products of two graphs.
\begin{theorem}\label{Cartesian.Prod.j=n}
Let $G$ realize $S_{\{i,j\}_{n}^{m}}$ for $n\geqslant 9$. Then $G$ is not the Cartesian product of two graphs.
\end{theorem}
\begin{proof}
Let $G$ realize $S_{\{i,j\}_{n}^{m}}$ with $n\geqslant9$, and suppose that $G=G_{1}\times G_{2}$. By Theorem~\ref{Cartesian.eigen},
all the eigenvalues of the graphs~$G_{1}$ and~$G_{2}$ belong to the Laplacian spectrum of~$G$, so $G_{1}$ or $G_{2}$ must have only integer eigenvalues.
Moreover, both graphs must have simple eigenvalues. Otherwise, $G$ has at least two multiple eigenvalues according to~\eqref{Cartesian.Prod.Spectrum}.

Thus, the graphs  $G_{1}$ and $G_{2}$  have simple integer Laplacian eigenvalues, so that they realize some sets $S_{l,r}$ and $S_{q,p}$, respectively,
and $n=rp\geqslant9$. By~\eqref{Cartesian.Prod.Spectrum} the spectrum of $G$ contains the eigenvalues $(l-1)+(q+1)=q+l$ and $(l+1)+(q-1)=q+l$, so
the Laplacian eigenvalue~$q+l$ of the graph $G$ has multiplicity at least~$2$. Let us show that the graph $G$ has at least one more multiple eigenvalue.

Indeed, if $r=p$, then the graph $G$ has a multiple eigenvalue $r$ because $0+r$ and $r+0$ are both the Laplacian eigenvalues of~$G$
by~\eqref{Cartesian.Prod.Spectrum}. If $r\neq p$ (say, $r<p$) and $r\in S_{q,p}$, then the graph $G$ again has two equal Laplacian
eigenvalues $0+r$ and $r+0$.

If $r\neq p$ (say, $r<p$) and $r\notin S_{q,p}$, then $r=q$ and $r-1\in S_{q,p}$. Now if $l\neq r-1$,
then by Theorem~\ref{Cartesian.eigen}, the number $r-1$ is a Laplacian eigenvalue of~$G$ of multiplicity at least $2$. If $l=r-1$
and $r>2$, then $r-2$ belong to $S_{l,r}$ and $S_{q,p}$, so $r-2+0$ and $0+r-2$ are both belong to the Laplacian spectrum of~$G$.
Finally, if $l=r-1$ and $r=2$, we have that $p\geqslant 5$ as $G$ is of order $rp\geqslant9$. Therefore, the eigenvalues $3$ and
$5$ are in the Laplacian spectrum of the graph $G_{2}$ realizing $S_{2,p}$ in this case. According to Theorem~\ref{Cartesian.eigen},
the graph $G$ has a multiple Laplacian eigenvalue $5$, since $0+5=2+3=5$ belongs to the Laplacian spectrum of~$G$.
	
Thus, if a graph $G$ realizing $S_{\{i,j\}_{n}^{m}}$ is a Cartesian product of two graphs, then it must have at least two multiple
Laplacian eigenvalues.
\end{proof}
\begin{remark}\label{Remark.Cartesian.product}
For $n\leqslant 8$, the only sets $S_{\{i,j\}_{n}^{m}}$ realized by Cartesian products are $S_{\{1,3\}_{4}^{2}}$,   $S_{\{4,6\}_{6}^{3}}$, and $S_{\{7,8\}_8^3}$.
For $n=4,6$ it follows from Tables~\ref{Table.2} and~\ref{Table.4} in Appendix~\ref{Tables.constructions}. For $n=8$, we have
only $5$ Cartesian products $K_2$ with a Laplacian integral graph on $4$ vertices. Among them there is only $S_{\{7,8\}_8^3}$ is of
kind $S_{\{i,n\}_n^m}$.
\end{remark}

\setcounter{equation}{0}
\section{Graphs realizing the sets $S_{\{i,j\}_{n}^{n}}$}\label{section.results.1}

In this section, we describe the graphs realizing the sets $S_{\{i,j\}_{n}^{n}}$, and first we note that if $G$ is a connected graph of order $n$ realizing a set $S_{j,n}$,
then according to Theorem~\ref{Join.Thm}, the graph $K_{1}\vee G$ realizes the set
\begin{equation*}
S_{\{1,j+1 \}_{n+1}^{n+1}}=\{ 0,2,3,\ldots,j-1,j,j+2,\ldots,n,n+1,n+1\}.
\end{equation*}
Thus, we come to the following lemma (cf.~\cite[Lemma~2.5]{FallatKirkland_et_al_2005}).
\begin{lemma}\label{Join.lemma}
If the set $S_{j,n}$ is Laplacian realizable, then so is $S_{\{1,j+1\}_{n+1}^{n+1}}$.
\end{lemma}

Now we are in a position to describe all the Laplacian realizable sets $S_{\{i,j\}_n^n}$ (cf.~Proposition~\ref{All.Sin}).
\begin{theorem}\label{Case_m=n}
Suppose $n\geqslant 3$. The only Laplacian realizable sets $S_{\{i,j\}_{n}^{n}}$, $i<j$ are the following.
\begin{itemize}
\item[(i)] If $n\equiv 0\mod 4$, then for each $k=1,2,\ldots,\frac{n-2}{2}$, $S_{\{1,2k+1\}_{n}^{n}}$ is Laplacian realizable;
\item[(ii)] If $n\equiv 1\mod 4$, then for each  $k=1,2,\ldots,\frac{n-3}{2}$, $S_{\{1,2k+1\}_{n}^{n}}$ is Laplacian realizable;
\item[(iii)] If $n\equiv 2\mod 4$, then for each $ k=1,2,\ldots,\frac{n-2}{2}$, $S_{\{1,2k\}_{n}^{n}}$ is Laplacian realizable;
\item[(iv)] If $n\equiv 3\mod 4$, then for each $ k=1,2,\ldots,\frac{n-1}{2}$, $S_{\{1,2k\}_{n}^{n}}$ is Laplacian realizable.
\end{itemize}
\end{theorem}
\begin{proof}
By Theorem~\ref{corol.Join.n.n}, if $S_{\{i,j\}_{n}^{n}}$ is Laplacian realizable then $i=1$, so the sets $S_{\{i,j\}_{n}^{n}}$ are not Laplacian realizable for $i>1$.
In~\cite[p.~286--289]{CvetkovicRowlinson_2010} (see also Tables~\ref{Table.2}--\ref{Table.4}) the authors found the Laplacian spectra of all graphs up to order~$5$. From that paper it follows that for $n\leqslant5$ the only
Laplacian realizable $S_{\{i,j\}_{n}^{n}}$ sets are $S_{\{1,2\}_{3}^{3}}$, $S_{\{1,3\}_{4}^{4}}$ and~$S_{\{1,3\}_{5}^{5}}$ covered
by items~(iv),(i) and (ii), respectively. Suppose now that $n\geqslant 6$.
\begin{itemize}							
\item[(i)] If $n\equiv 0\mod 4$, then $n-1\equiv3\mod4$. Therefore, by Proposition~\ref{All.Sin}~(iv), the set $S_{2k,n-1}$ is Laplacian realizable for $k=1,2,\ldots,\frac{n-2}{2}$.
Now from Lemma~\ref{Join.lemma}, we have that $S_{\{1,2k+1\}_{n}^{n}}$ is Laplacian realizable for any $k=1,2,\ldots,\frac{n-2}{2}$. At the same time,
since $n$ is even and $m=n$, $S_{\{1,2k\}_{n}^{n}}$ are not Laplacian realizable by Theorem~\ref{Condition.m}~(i).

\item[(ii)]  If $n\equiv 1\mod 4$, then $n-1\equiv0\mod4$, so for each $k=1,2,\ldots,\frac{n-3}{2}$, $S_{2k,n-1}$ is Laplacian realizable by Proposition~\ref{All.Sin}~(i).
Lemma~\ref{Join.lemma} implies that $S_{\{1,2k+1\}_{n}^{n}}$ is Laplacian realizable for any $k=1,2,\ldots,\frac{n-3}{2}$, while Theorem~\ref{Condition.m}~(ii)
gives that $S_{\{1,2k\}_{n}^{n}}$ is not Laplacian realizable, since $n$ is odd, so is $m=n$.
%
%
%
\end{itemize}	
The cases (iii) and (iv) can be proved analogously.
\end{proof}

Our next result discusses the structure of graphs realizing $S_{\{1,j\}_{n}^{n}}$ for various possible~$j$.
\begin{theorem}\label{Constru.Double.Eig}
Let $G$ be a graph of order $n$, $n\geqslant5$.
\begin{itemize}
\item [(a)]  The graph $G$ realizes $S_{\{1,2\}_{n}^{n}}$ if and only if $G$ is formed in one of the following two ways:
\begin{itemize}
\item[(i)]  $G=P_3\vee\left(K_{1}\cup H\right)$, where $H$ realizes $S_{n-5,n-4}$ and $P_{3}$ is the path graph on $3$ vertices;
\item[(ii)]  $G=K_2\vee H$, where $H$ realizes $S_{n-2,n-2}$.
\end{itemize}
\item[(b)] If $3\leqslant j\leqslant n-2$, then $G$ realizes $S_{\{1,j\}_{n}^{n}}$ if and only if $G=K_{2}\vee\left(K_{1}\cup H\right)$, where the graph $H$ realizes~$S_{j-2,n-3}$.
\item [(c)]  The graph $G$ realizes $S_{\{1,n-1\}_{n}^{n}}$ if and only if $G$ is formed in one of the following two ways:
\begin{itemize}
\item[(i)]  $G=K_2\vee(K_{2}\cup H)$, where $H$ realizes $S_{2,n-4}$;
\item[(ii)]  $G=K_2\vee(K_{1}\cup H)$, where $H$ realizes $S_{n-3,n-3}$.
\end{itemize}
\end{itemize}
\end{theorem}
\begin{proof}

For $n=5$ there is only one graph realizing a set  $S_{\{1,j\}_{n}^{n}}$, according to Table~\ref{Table.3}. This graph
is $G=K_2\vee(K_1\cup K_2)=K_2\vee\overline{P_3}$ whose Laplacian spectrum
is $S_{\{1,3\}_{5}^{5}}$. This graph satisfies condition $(b)$ of the theorem.

For $n=6$ there are only two graphs realizing sets $S_{\{1,j\}_{n}^{n}}$ according to Table~\ref{Table.4}.
The graph $G_1=P_3\vee(K_1\cup K_2)=P_3\vee\overline{P_3}$ has the Laplacian spectrum $S_{\{1,2\}_6^6}$ and satisfies
condition (a)(i), while the graph $G_2=K_{2}\vee\left(K_{1}\cup P_3\right)$ with Laplacian spectrum
$S_{\{1,4\}_6^6}$ satisfies condition (b). Thus, the theorem is true for $n=5$ and $6$.

\vspace{2mm}

Let now $n\geqslant7$, and suppose that $G$ realizes $S_{\{1,2\}_{n}^{n}}$. Then by Theorem~\ref{Complement.spect},
one has $\sigma_L(\overline{G})= \{0\}\cup\{0\}\cup S_{n-2,n-2}$. Therefore, one can represent the complement of the
graph $G$ as follows: $\overline{G}=K_1\cup\overline{F}$, where $\overline{F}$ is disconnected and
$\sigma_L(\overline{F})=\{0\}\cup S_{n-2,n-2}$. Again by Theorem~\ref{Complement.spect}, we obtain
$\sigma_L(F)=S_{1,n-1}$. According to Proposition~\ref{Propos:Constr.S_i.n.222}, there are only two possibilities
to construct the graph $F$.

If $F=(K_{1}\cup K_{1})\vee (K_{1}\cup H)$, where $H$ is a graph on $n-4$ vertices realizing $S_{n-5,n-4}$
(see Proposition~\ref{Propos:Constr.S_i.n.222}~(i)), then we obtain $\overline{G}=K_1\cup\overline{K_{1}\cup K_{1}}\cup\overline{K_{1}\cup H}$.
It is clear that $K_1\cup\overline{K_{1}\cup K_{1}}=K_1\cup K_2=\overline{P_3}$, where $P_3$ is the path graph
on $3$ vertices. Thus, we obtain that $G=P_3\vee(K_{1}\cup H)$ where $H$ realizes the set $S_{n-5,n-4}$.

If $F=K_{1}\vee H$, where  $H$ is a graph on $n-2$ vertices that realizes $S_{n-2,n-2}$ (see
Proposition~\ref{Propos:Constr.S_i.n.222}~(ii)), then $\overline{G}=K_1\cup\overline{K_{1}\vee H}=(K_1\cup K_1)\cup \overline{H}$.
Therefore, $G=K_2\vee H$, since $K_2=\overline{K_1\cup K_1}$. It can be easily checked that the graphs
mentioned in conditions $(a)(i)$ and $(a)(ii)$ realize the set $S_{\{1,2\}_{n}^{n}}$, so the statement
$(a)$ of the theorem is proved.

\vspace{2mm}

Let now $n\geqslant7$ and $3\leqslant j\leqslant n-2$. If $G$ realizes $S_{\{1,j\}_{n}^{n}}$, then from Theorem~\ref{Complement.spect}, we obtain $\sigma_L{(\overline G)}=\{0\}\cup\{0\}\cup S_{n-j,n-2}$.
So the complement of the graph $G$ can be represented as follows: $\overline{G}=K_1\cup\overline{F}$, where $\overline{F}$ is disconnected and $\sigma_L(\overline{F})=\{0\}\cup S_{n-j,n-2}$.
By Theorem~\ref{Complement.spect}, we obtain $\sigma_L(F)=S_{j-1,n-1}$. According to Proposition~\ref{Propos:Constr.S_i.n.223}, we have
$F=K_{1}\vee (K_{1}\cup H)$, where $H$ is a graph on $n-3$ vertices that realizes $S_{j-2,n-3}$. Therefore, $\overline{G}=K_1\cup K_{1}\cup\overline{K_{1}\cup H}
=\overline{K_2}\cup\overline{K_{1}\cup H}$. Thus, $G=K_2\vee(K_{1}\cup H)$ where $H$ is a graph realizing $S_{j-2,n-3}$. Converse statement $(b)$ can be easily checked.

\vspace{2mm}

Finally, let $n\geqslant7$, and suppose that $G$ realizes $S_{\{1,n-1\}_{n}^{n}}$. Then by Theorem~\ref{Complement.spect}, one has $\sigma_L(\overline{G})= \{0\}\cup\{0\}\cup S_{1,n-2}$. Therefore, one can represent
the complement of the graph $G$ as follows: $\overline{G}=K_1\cup\overline{F}$, where $\overline{F}$ is disconnected and $\sigma_L(\overline{F})=\{0\}\cup S_{1,n-2}$.
Again from Theorem~\ref{Complement.spect}, we have $\sigma_L(F)=S_{n-2,n-1}$. According to Proposition~\ref{Propos:Constr.S_i.n.224}, there are only two possibilities
to construct the graph $F$.

If $F=K_{1}\vee (K_{2}\cup H)$, where $H$ is a graph on $n-4$ vertices realizing $S_{2,n-4}$ (see Proposition~\ref{Propos:Constr.S_i.n.224} (i)),
then we obtain $\overline{G}=K_1\cup K_{1}\cup\overline{K_{2}\cup H}$. Thus, we obtain that $G=K_2\vee(K_{2}\cup H)$ where $H$ realizes the set $S_{2,n-4}$.

If $F=K_{1}\vee (K_{1}\cup H)$, where  $H$ is a graph on $n-3$ vertices that realizes $S_{n-3,n-3}$ (see Proposition~\ref{Propos:Constr.S_i.n.224} (ii)), then
$\overline{G}=K_1\cup\overline{K_{1}\vee (K_{1}\cup H)}=(K_1\cup K_1)\cup \overline{K_{1}\cup H}$. Therefore, $G=K_2\vee(K_{1}\cup H)$. It can be easily checked
that the graphs mention in conditions $(c)(i)$ and $(c)(ii)$ realize the set $S_{\{1,n-1\}_{n}^{n}}$, so the statement $(c)$ of the theorem is proved.
\end{proof}

\begin{figure}[h]
\centering
\includegraphics[width=7cm,height=3cm]{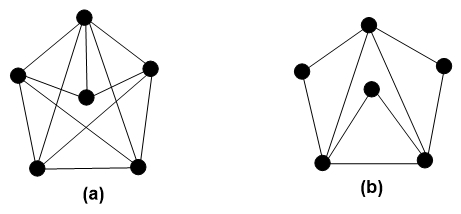}
\caption{Graphs realizing $S_{\{1,2\}_6^6}$ and $S_{\{1,4\}_6^6}$ respectively.}\label{Isomo.Graph}
\end{figure}

\begin{remark}
Note that the only graphs of orders 3 and 4 realizing sets $S_{\{1,j\}_n^n}$ are $K_3$ and $K_2\vee2K_1$. They satisfy the conditions of
Theorem~\ref{Constru.Double.Eig}, as well, if we formally set that the isolated vertex has the Laplacian spectrum $S_{1,1}:=\{0\}$.
\end{remark}
Figure~\ref{Isomo.Graph} illustrates the graphs realizing sets $S_{\{1,j\}_n^n}$ for $n=6$.

In Theorem~\ref{Case_m=n} we have completely resolved the existence of graphs realizing
the spectrum $S_{\{i,j\}_n^n}$, where $1\leqslant i<j<n$. Furthermore, it is easily deduced from
Theorem~\ref{Constru.Double.Eig} that if the sets $S_{n,n}$ were not realizable for any $n$, then there is a
unique graph, which realizes $S_{\{i,j\}_n^n}$, for $1\leqslant i<j<n$. We remind the reader that
according to the so-called $S_{n,n}$-conjecture, the sets $S_{n,n}$ are not realizable,
see~\cite{FallatKirkland_et_al_2005,GoldbergerNeumann_2013}.

\setcounter{equation}{0}

\section{Graphs realizing the sets $S_{\{i,j\}_{n}^{n-1}}$ }\label{section.results.2}

In this section, we provide a criterion for the sets $S_{\{i,j\}_{n}^{n-1}}$ to be Laplacian
realizable and describe realizing these sets.

We start from a necessary condition for $S_{\{i,j\}_{n}^{n-1}}$ to be Laplacian realizable.
\begin{theorem}\label{condition.n-1}
If $S_{\{i,j\}_{n}^{n-1}}$ is Laplacian realizable then the number $i$ is either $1$ or $2$.
\end{theorem}
Thus, the theorem claims that the sets $S_{\{i,j\}_{n}^{n-1}}$ are not Laplacian realizable
for $i\geqslant 3$.
\begin{proof}
Suppose, on the contrary, that $S_{\{i,j\}_{n}^{n-1}}$ is Laplacian realizable for some
$i\geqslant3$ (here $j>i$), and $G$ is a graph realizing this set. Then by Theorem~\ref{join.spectrum.1}, $G$ is
a join of two graphs, say, $G=F\vee J$.
Let the order of~$F$ be $p$ and the order of $J$ be $n-p$. We denote the Laplacian spectra
$F$ and $J$ as follows
\begin{equation*}
0=\mu_{1}\leqslant\mu_{2}\leqslant\mu_{3}\ldots\leqslant \mu_{p-1}\leqslant\mu_{p}\quad\text{and}\quad 0=\lambda_{1}\leqslant\lambda_{2}\leqslant\lambda_{3}\ldots\leqslant\lambda_{n-p-1}\leqslant\lambda_{n-p}.
\end{equation*}
Then by Theorem~\ref{Join.Thm}, the Laplacian spectrum of~$G$ has the form
\begin{equation}\label{spectrum.join.2}
\sigma_{L}(G)=\{0,(n-p)+\mu_{2},\ldots,(n-p)+\mu_{p},p+\lambda_{2},\ldots,p+\lambda_{n-p-1},p+\lambda_{n-p},n\},
\end{equation}
where the eigenvalues are not in the increasing order.

Since $1$ is a simple Laplacian eigenvalue of~$G$, one has either $(n-p)+\mu_{2}=1$ or $p+\lambda_{2}=1$. Without loss of generality,
we can suppose that $(n-p)+\mu_{2}=1$ and $p+\lambda_{2}\neq 1$. This implies that $n-p=1$ and $\mu_{2}=0$, so the graph $J$
is an isolated vertex with Laplacian spectrum $\sigma_L(J)=\{0\}$ while the graph $F$ is a union of two disjoint connected graphs,
$F=H_1\cup H_2$. Now we recall that  $n-1$ is the only double Laplacian eigenvalue of~$G$ by assumption, so from~\eqref{spectrum.join.2},
we have $n-p+\mu_{p-1}=n-p+\mu_{p}=n-1$ that implies $\mu_{p-1}=\mu_{p}=p-1=n-2$. Thus, from Proposition~\ref{Prop.max.eig} we obtain that
$H_1=K_1$, and the graph $H_2$ is of order~$n-2$ with Laplacian eigenvalue $n-2$ of multiplicity $2$.

On the other hand, the number $2$ is also a simple Laplacian eigenvalue of $G$, so $(n-p)+\mu_{3}=2$ and hence $\mu_{3}=1$. Thus, $H_2$
has a Laplacian eigenvalue $1$ that contradicts Theorem~\ref{corol.Join.n.n}. Therefore, the sets~$S_{\{i,j\}_{n}^{n-1}}$ are
not Laplacian realizable for $i\geqslant3$.	
\end{proof}				
%

\subsection{Graphs realizing the sets $S_{\{2,j\}_{n}^{n-1}}$}\label{subsection1:results.2}

First we describe the sets $S_{\{2,j\}^{n-1}_n}$, because they are used in the description of the sets $S_{\{1,j\}^{n-1}_n}$. The following
theorem lists all the Laplacian realizable sets $S_{\{2,j\}^{n-1}_n}$.
\begin{theorem}\label{double.m=n-1}
Let $n\geqslant 3$, and let $G$ be a connected graph of order $n$. The only Laplacian realizable sets~$S_{\{2,j\}_n^{n-1}}$ are the following:
\begin{itemize}
\item[(i)] If $n\equiv 0\mod 4$, then for each  $k=1,2,\ldots,\dfrac{n-4}{2}$, $S_{\{2,2k+1\}_n^{n-1}}$ is Laplacian realizable;
\item[(ii)] If $n\equiv 1\mod 4$, then for each $k=1,2,\ldots,\dfrac{n-3}{2}$, $S_{\{2,2k+1\}_n^{n-1}}$ is Laplacian realizable;
\item[(iii)] If $n\equiv 2\mod 4$, then for each $ k=1,2,\ldots,\dfrac{n-4}{2}$,
$S_{\{2,2k+2\}_n^{n-1}}$  is Laplacian realizable;
\item[(iv)] If $n\equiv 3\mod 4$,then for each $ k=1,2,\ldots,\dfrac{n-5}{2}$, $S_{\{2,2k+2\}_n^{n-1}}$ is Laplacian realizable.
\end{itemize}
\end{theorem}
\begin{proof}
For $n=5$, the only Laplacian realizable $S_{\{i,j\}_{n}^{n-1}}$ set is $S_{\{2,3\}_5^4}$, see Table~\ref{Table.3}
in Appendix~\ref{Tables.constructions}. This agrees with the condition (ii) of the theorem. Suppose now that $n\geqslant 6$.
\begin{itemize}
\item[(i)] If $n\equiv 0 \mod 4$, then $n-2\equiv 2\mod 4$, then for $k=1,2,\ldots,\dfrac{n-4}{2}$, the sets $S_{\{1,2k\}_{n-2}^{n-2}}$
are Laplacian realizable by Theorem~\ref{Case_m=n}~(iii). From Lemma~\ref{join.union} it follows that $S_{\{2,2k+1\}_n^{n-1}}$
is Laplacian realizable for each $k=1,2,\ldots,\dfrac{n-4}{2}$. At the same time, the sets $S_{\{2,2k\}_n^{n-1}}$ are not Laplacian
realizable for any $k$ by Theorem~\ref{Condition.m}~(i), since the double eigenvalue $m=n-1$ is an odd number in this case.
\item[(ii)] For $n\equiv 1 \mod 4$, we have $n-2\equiv 3\mod 4$. Thus,  for $k=1,2,\ldots,\dfrac{n-3}{2}$, the sets $S_{\{1,2k\}_{n-2}^{n-2}}$
are Laplacian realizable by Theorem~\ref{Case_m=n}~(iv). Now Lemma~\ref{join.union} implies the $S_{\{2,2k+1\}_n^{n-1}}$ are
Laplacian realizable for $k=1,2,\ldots,\dfrac{n-3}{2}$. As $m=n-1$ is even, from Theorem~\ref{Condition.m}~(ii) it follows that
the sets $S_{\{2,2k\}_n^{n-1}}$  are not Laplacian realizable for any $k$.
\end{itemize}

The cases (iii) and (iv) can be proved analogously with use of Theorems~\ref{Condition.m} and~\ref{Case_m=n} and Lemma~\ref{join.union}.

\end{proof}
In the following theorem, we discuss construction of graphs realizing the sets $S_{\{2,j\}_n^{n-1}}$.
\begin{theorem}\label{Thm.case.2.m=n-1}
Let $n\geqslant 5$, and let $G$ be a connected graph of order $n$. Then $G$ realizes
$S_{\{2,j\}_n^{n-1}}$ if and only if $G=K_{1}\vee\left(K_{1}\cup H\right)$,
where $H$ is a graph on $n-2$ vertices realizing $S_{\{1,j-1\}_{n-2}^{n-2}}$;
\end{theorem}
\begin{proof}
Let $G$ realize $S_{\{2,j\}_n^{n-1}}$. Then $G=G_1\vee G_2$ by Theorem~\ref{Thm.Join.n},
so $\overline{G}=\overline{G_{1}}\cup\overline{G_{2}}$. From Theorem~\ref{Complement.spect} it follows that
$\sigma_L(\overline G)=\{0\}\cup S_{\{n-j,n-2\}_{n-1}^1}$. Thus, by Proposition~\ref{Prop.max.eig} we obtain that $G_1=K_1$,
and $G_2$ is of order $n-1$, so that $\sigma_L(\overline{G_2})=S_{\{n-j,n-2\}_{n-1}^1}$. Using Theorem~\ref{Complement.spect},
we get $\sigma_L(G_2)=\{0\}\cup S_{\{1,j-1\}_{n-2}^{n-2}}$. Again Proposition~\ref{Prop.max.eig}
gives us that $G_2=K_1\cup H$, where $H$ is a graph on $n-2$ vertices realizing~$S_{\{1,j-1\}_{n-2}^{n-2}}$.
Consequently, $G=K_1\vee \left(K_{1}\cup H\right)$, as required.

Conversely, if $G=K_{1}\vee(K_{1}\cup H)$, where $H$ is a graph on $n-2$ vertices realizing~$S_{\{1,j-1\}_{n-2}^{n-2}}$,
then from Theorem~\ref{Join.Thm} it follows that $G$ realizes $S_{\{2,j\}_n^{n-1}}$.
\end{proof}

\begin{remark}
Theorems~\ref{double.m=n-1} and~\ref{Thm.case.2.m=n-1} completely resolve the existence of graphs realizing
the spectrum~$S_{\{2,j\}_n^n}$, where $3\leqslant j<n$. The realizability of the sets $S_{\{2,n\}_n^{n-1}}$
and, in general, the sets $S_{\{i,n\}_n^{m}}$ is discussed in Section~\ref{section.conj}.
\end{remark}

\subsection{Graphs realizing the sets $S_{\{1,j\}_{n}^{n-1}}$}\label{subsection2:results.2}

Now we are in position to study the realizability of sets~$S_{\{1,j\}_n^{n-1}}$.
It turns out that for a fixed $n$, there is only one Laplacian realizable set of this kind.

\begin{theorem}\label{Theorem.S.1.j.n-1.n}
Let $G$ be a simple connected graph of order $n$, $n\geqslant6$.
\begin{itemize}
\item[(i)] For $n\equiv 0$ or $1\mod 4$, the set $S_{\{1,j\}_n^{n-1}}$ is Laplacian realizable if and only if $j=2$.
\item[(ii)] For $n\equiv 2$ or $3\mod 4$, the set $S_{\{1,j\}_n^{n-1}}$ is Laplacian realizable if and only if $j=3$.
\end{itemize}
\end{theorem}
\begin{proof}

\noindent (i) If $n\equiv 0$ or $1\mod 4$ (so that $n\geqslant8$), then from Theorem~\ref{Condition.m} it follows that the sets $S_{\{1,2k+1\}_n^{n-1}}$ are
not Laplacian realizable for any $k$. Let us now find out which sets of the form $S_{\{1,2k\}_n^{n-1}}$ are realizable.

Suppose first that $k\geqslant2$. Let a graph $G$ realize $S_{\{1,2k\}_n^{n-1}}$. Then we have
\begin{equation}\label{Theorem.S.1.j.n-1.n.proof.1}
\sigma_L(\overline{G})=\{0,0,1,1,2,\ldots,n-(2k+1),n-(2k-1),\ldots,n-3,n-2\}.
\end{equation}
Thus, the graph $\overline{G}$ is a union of two disjoint connected graphs, $\overline{G}=\overline{H_1}\cup \overline{H_2}$. From
Proposition~\ref{Prop.max.eig} it follows that one of the components, say, $H_2$ must be of order at least $n-2$.

If $H_{2}$ has $n-2$ vertices, then $H_1$ has two vertices, and therefore $\sigma_L(\overline{H_{1}})=\{0,2\}$. By
Theorem~\ref{Complement.spect}, we obtain $\sigma_L(H_{1})=\{0,0\}$, so $H_1=2K_1$. Now from~\eqref{Theorem.S.1.j.n-1.n.proof.1}
we obtain $\sigma_L(\overline{H_2})=S_{\{2,n-2k\}_{n-2}^{1}}$, and Theorem~\ref{Complement.spect} implies
$\sigma_L(H_2) =\{0,0,1,2,\ldots,2k-4,2k-3,2k-1,2k,\ldots,n-5,n-3,n-3\}$. By Proposition~\ref{Prop.max.eig}
we have that $H_2=K_1\cup H_3$ where $\sigma_L(H_3)=S_{\{2k-2,n-4\}_{n-3}^{n-3}}$, so $1\in \sigma_L(H_3)$.
This contradicts Theorem~\ref{corol.Join.n.n}.
	
If $H_{2}$ has $n-1$ vertices, then $H_1=K_1$, so that $\overline{G}=K_1\cup \overline{H_2}$ with $\sigma_L(\overline{H_2})=S_{\{n-2k,n-1\}_{n-1}^{1}}$.
Now Theorem~\ref{Complement.spect} gives us $\sigma_L(H_2)=S_{\{2k-1,n-1\}_{n-1}^{n-2}}$. However,
for $k\geqslant2$ such sets are not Laplacian realizable according to Theorem~\ref{condition.n-1}. Thus,
we obtained that the sets $S_{\{1,2k\}_n^{n-1}}$ are not Laplacian realizable for~$k\geqslant2$.
	
Assume now that the graph $G$ realizes $S_{\{1,2\}_n^{n-1}}$.
According to Theorem~\ref{Complement.spect},
\begin{equation*}
\sigma_L(\overline{G})= \{0,0,1,1,2,\ldots,n-4,n-3\}.
\end{equation*}
Thus, the graph $\overline{G}$ is a union of two disjoint connected graphs, $\overline{G}=\overline{H_1}\cup\overline{H_2}$. From
Proposition~\ref{Prop.max.eig} it follows that one of the components, say, $H_2$ must be of order at least $n-3$.

Let $H_{2}$ have $n-3$ vertices. Then $H_{1}$ has three vertices and the only possible spectrum\footnote{The
only connected graphs on $3$ vertices are $P_{3}$ and $K_{3}$, but the Laplacian spectrum of $K_{3}$ contains double $3$.}
of $\overline{H_1}$ is $\{0,1,3\}$ that implies $H_1=K_1\cup K_2$ by Theorem~\ref{Complement.spect}. Therefore, $\sigma_L(\overline{H_2})=S_{3,n-3}$. By Theorem~\ref{Complement.spect},
one obtains  $\sigma_L(H_2)=\{0\}\cup S_{n-6,n-4}$. By Theorem~\ref{All.Sin} these sets are Laplacian realizable
for $n\equiv 0$ or $1\mod 4$. Thus, the graph $G=(K_1\cup K_2)\vee(K_{1}\cup H)$, where the graph $H$ realizes $S_{n-6,n-4}$.
	
If $H_{2}$ has $n-2$ vertices, then $H_{1}$ has two vertices and $\sigma_L(\overline{H_{1}})=\{0,2\}$, so that $H_1=2K_1$.
Hence $\sigma_L(\overline{H_{2}})=S_{\{2,n-2\}_{n-2}^{1}}$. By Theorem~\ref{Complement.spect}, we have that
the graph $H_2$ must realize the set $S_{\{n-4,n-2\}_{n-2}^{n-3}}$. However, this set
is not Laplacian realizable for $n\geqslant8$ according to Theorem~\ref{condition.n-1}. Consequently, $H_2$ cannot be of
order $n-2$.
	
Suppose now that $H_{2}$ has $n-1$ vertices. Then $H_{1}=K_1$, and $\sigma_L(\overline{H_2})=S_{\{n-2,n-1\}_{n-1}^{1}}$.
Theorem~\ref{Complement.spect} gives $\sigma_L(H_2)=S_{\{1,n-1\}_{n-1}^{n-2}}$. However, this set
is not Laplacian realizable as we established above.
	
Thus, we obtained that if $S_{\{1,j\}_n^{n-1}}$ is Laplacian realizable then $j=2$, and the graph
$G=(K_1\cup K_2)\vee(K_1\cup H)$ with $\sigma_L(H)=S_{n-6,n-4}$ realizes $S_{\{1,2\}_n^{n-1}}$.
Converse statement is obvious.

\vspace{2mm}	

\noindent (ii) Let now $n\equiv 2$ or $3\mod 4$. Then from Theorem~\ref{Condition.m} it follows that the
sets $S_{\{1,2k\}_n^{n-1}}$ are not Laplacian realizable for any $k$. We now find out which
sets of the form $S_{\{1,2k+1\}_n^{n-1}}$ are realizable.

As above we first consider the case $k\geqslant2$. Let the graph $G$ realize $S_{\{1,2k+1\}_n^{n-1}}$.
Then we have
\begin{equation*}
\sigma_L(\overline{G})=\{0,0,1,1,2,\ldots,n-2k-1,n-2k+1,\ldots,n-3,n-2\}.
\end{equation*}
Thus, the graph $\overline{G}$ is a union of two disjoint connected graphs, $\overline{G}=\overline{H_1}\cup \overline{H_2}$. From
Proposition~\ref{Prop.max.eig} it follows that one of the components, say, $H_2$ must be of order at least $n-2$.

If $H_2$ is of order $n-2$, then the Laplacian spectrum of $\overline{H_1}$ must be $\{0,2\}$, so that $H_1=2K_1$.
Consequently, $\sigma_L(\overline{H_2})=S_{\{2,n-2k-1\}_{n-2}^{1}}$, and by Theorem~\ref{Complement.spect}
$\sigma_L(H_2)=\{0\}\cup S_{\{2k-1,n-4\}_{n-3}^{n-3}}$. By Proposition~\ref{Prop.max.eig} we have $H_2=K_1\cup H_3$
where $\sigma_L(H_3)=S_{\{2k-1,n-4\}_{n-3}^{n-3}}$. Since the set $S_{\{2k-1,n-4\}_{n-3}^{n-3}}$ is not Laplacian
realizable for $k\geqslant2$ (Theorem~\ref{corol.Join.n.n}), $H_2$ cannot be of order $n-2$.

Let now the order of $H_2$ be $n-1$. Then we have $H_1=K_1$, and $\sigma_L(\overline{H_2})=S_{\{n-2k-1,n-1\}_{n-1}^{1}}$,
so that $\sigma_L(H_2)=S_{\{2k,n-1\}_{n-1}^{n-2}}$ by Theorem~\ref{Complement.spect}. Now from Theorem~\ref{condition.n-1}
it follows that for $k\geqslant2$ such sets are not realizable. Thus, if $S_{\{1,2k+1\}_n^{n-1}}$ is Laplacian realizable,
then $k$ can only be equal $1$ (if any).

Suppose now that a graph $G$ realizes the set $S_{\{1,3\}_n^{n-1}}$ for $n\equiv 2$ or $3\mod 4$. Then its complement
has the following Laplacian spectrum:
\begin{equation*}
\sigma_L(\overline{G})=\{0,0,1,1,2,\ldots,n-4,n-2\}.
\end{equation*}
As above we obtain that $\overline{G}=\overline{H_1}\cup \overline{H_2}$ where one of the components, say, $H_2$ must be of order at least $n-2$.

If the order of $H_2$ equals $n-2$, then necessarily $\sigma_L(\overline{H_1})=\{0,2\}$, and $H_1=2K_1$ while $\sigma_L(\overline{H_2})=S_{\{2,n-3\}_{n-2}^{1}}$.
Thus, by Theorem~\ref{Complement.spect} we have $\sigma_L(H_2)=\{0\}\cup S_{\{1,n-4\}_{n-3}^{n-3}}$. By Proposition~\ref{Prop.max.eig} we have $H_2=K_1\cup H_3$
where $\sigma_L(H_3)=S_{\{1,n-4\}_{n-3}^{n-3}}$. Finally, we obtain $G=2K_1\vee(K_1\cup F)$ where $\sigma_L(F)=S_{\{1,n-4\}_{n-3}^{n-3}}$.

Let $H_2$ be of order $n-1$. Then $H_1=K_1$, and $\sigma_L(\overline{H_2})=S_{\{n-3,n-1\}_{n-1}^{1}}$, so that $\sigma_L(H_2)=S_{\{2,n-1\}_{n-1}^{n-2}}$.
This case is realizable if the set $S_{\{2,n-1\}_{n-1}^{n-2}}$ is Laplacian realizable.
\end{proof}

From the proof of Theorem~\ref{Theorem.S.1.j.n-1.n} one can easily see the structure of graphs realizing~$S_{\{1,j\}_{n}^{n-1}}$.
\begin{theorem}\label{Theorem.S.1.j.n-1.n_constr}
Let $G$ be a graph of order $n$, $n\geqslant6$.
\begin{itemize}
\item [(a)]  The graph $G$ realizes $S_{\{1,2\}_{n}^{n-1}}$ if and only if $n\equiv 0$ or $1\mod 4$, and
$G=(K_1\cup K_2)\vee(K_1\cup H)$, where~$H$ realizes $S_{n-6,n-4}$;
\item [(b)] The graph $G$ realizes $S_{\{1,3\}_{n}^{n-1}}$ if and only if $n\equiv 2$ or $3\mod 4$, and
$G$ is formed in one of the following two ways:
\begin{itemize}
\item[(i)] $G=\left(K_{1}\cup K_{1}\right)\vee\left(K_{1}\cup H\right)$, where the graph $H$ realizes
$S_{\{1,n-4\}_{n-3}^{n-3}}$;
\item[(ii)] $G=K_{1}\vee F$, where the graph $F$ realizes $S_{\{2,n-1\}_{n-1}^{n-2}}$.
\end{itemize}
\end{itemize}
\end{theorem}
\begin{remark}
We have completely resolved the existence of graphs realizing
the spectrum $S_{\{1,j\}_{n}^{n-1}}$, where $j=1,3$ (other sets are not realizable).
Furthermore, it is easily deduced from Theorem~\ref{Theorem.S.1.j.n-1.n_constr} that
if the set~$S_{\{2,n\}_{n}^{n-1}}$ was not realizable for any $n$, then there is a
unique graph, which realizes $S_{\{1,3\}_{n}^{n-1}}$. We consider the realizability of
sets~$S_{\{i,n\}_n^m}$ in the next section.
\end{remark}

\setcounter{equation}{0}
\section{$S_{\{i,n\}_n^m}$-conjecture}\label{section.conj}

This section is devoting to describing necessary conditions on the set $S_{\{i,n\}_n^m}$ to
be Laplacian realizable. Since the results below do not resolve the question of the
realizability of $S_{\{i,n\}_n^m}$, we make Conjecture~\ref{Conjecture.S_nnm} at the
end of this section.

First, we note that by Theorem~\ref{Thm.Join.n}, any graph realizing $S_{\{i,n\}_n^m}$
is not a join of two graphs. However, it can be a Cartesian product of two graphs.
Indeed, the ladder graph on $6$ vertices, $P_{2}\times P_{3}$, realizes the set~$S_{\{4,6\}_6^3}$,
see Table~\ref{Table.4}. Moreover, from Theorem~\ref{Complement.spect} it immediately follows
that if $G$ realizes a~set~$S_{\{i,n\}_{n}^{m}}$, then its complement realizes a set of the same
kind, as well.
\begin{prop}\label{complement.S_i,n,n}
If $G$ realizes $S_{\{i,n\}_{n}^{m}}$, then $\overline{G}$ is a connected graph that realizes
$S_{\{n-i,n\}_{n}^{n-m}}$.
\end{prop}
Thus, together with $S_{\{4,6\}_6^3}$ we obtain that $S_{\{2,6\}_6^3}$ is Laplacian realizable,
see Figure~\ref{ladder1.pic} and Table~\ref{Table.4}.

Note that since the equalities $n-i=i$ and $n-m=m$ cannot hold simultaneously due to $i\neq m$, the graphs
realizing~$S_{\{i,n\}_{n}^{m}}$ are not self-complementary, so the number of graphs realizing
sets~$S_{\{i,n\}_{n}^{m}}$ is even (if finite).

According to \cite[p.~286--289]{CvetkovicRowlinson_2010} (see also Tables~\ref{Table.2}--\ref{Table.4} in
Appendix~\ref{Tables.constructions}), there are no graphs of order less than~$6$ realizing $S_{\{i,n\}_n^m}$
and exactly $2$ such graphs of order $6$. It can be proved that the sets $S_{\{i,n\}_n^m}$
are not Laplacian realizable for $n=7$ and for any prime $n$.
\begin{prop}\label{conj.facts}
If $n\geqslant7$ is a prime number, then $S_{\{i,n\}_n^m}$ is not Laplacian realizable.
\end{prop}
\begin{proof}
Let $G$ realizes a set $S_{\{i,n\}_n^m}$, and let the eigenvalues of $L(G)$ be given
by $0$, $\mu_2$, \ldots, $\mu_n$. According to the matrix tree
theorem, see, for instance,~\cite[p.~190]{CvetkovicRowlinson_2010}, the number of
spanning trees of $G$ is equal to
\begin{equation*}
\tau(G)=\dfrac{\lambda_2\cdots\lambda_n}{n}=\dfrac{(n-1)!\cdot m}{n\cdot i}.
\end{equation*}
As the number of spanning trees is an integer, $n$ must divide $\frac{(n-1)!\cdot m}{i}$ and so,
in particular, $n$ cannot be a prime number, a contradiction.
\end{proof}

From Theorem~\ref{Condition.m} we obtain another simple observation.

\begin{prop}\label{conj.facts.2}
Let $n$ be a positive integer, $n\geqslant6$. If $n\equiv 0$ or $1\mod 4$ and $i-m$ is odd,
then~$S_{\{i,n\}_n^m}$ is not Laplacian realizable.
\end{prop}

As we noticed above, graphs realizing $S_{\{i,n\}_n^m}$ are not joins but may be Cartesian products. However, by
Theorem~\ref{Cartesian.Prod.j=n} such graphs can be of order less than $9$. As we mentioned in Remark~\ref{Remark.Cartesian.product},
there are only~$3$ Cartesian products of order less than $9$, two of which realize sets of kind $S_{\{i,n\}_n^m}$. One of them
is $S_{\{4,6\}_6^3}$ mentioned above in this section. Another one is $S_{\{7,8\}_8^3}$. The graph realizing this set
is the Cartesian product of $K_2$ and $K_1\vee \overline{P_3}=K_1\vee(K_1\cup K_2)$ whose Laplacian spectrum is $S_{2,4}$.
By Theorem~\ref{complement.S_i,n,n}, the complement of $G=K_2\times[K_1\vee(K_1\cup K_2)]$ realizes a set of
kind~$S_{\{i,n\}_n^m}$ as well. Namely, $\sigma_L(\overline{G})=S_{\{1,8\}_8^5}$. Both graphs are depicted
in Figure~\ref{ladder2.pic}.

Thus, for $n\leqslant8$ we found $4$ Laplacian realizable sets $S_{\{i,n\}_n^m}$, two Cartesian products and their
complements. This fact together with Theorem~\ref{Cartesian.Prod.j=n} inspired us to make the following conjecture.
\begin{conjecture}\label{Conjecture.S_nnm}
For $n\geqslant4$, the only Laplacian realizable set of kind $S_{\{i,n\}_n^m}$ are $S_{\{4,6\}_6^3}$, $S_{\{2,6\}_6^3}$, $S_{\{7,8\}_8^3}$, and $S_{\{1,8\}_8^5}$.
\end{conjecture}
In light of this conjecture we have the following.
\begin{conjecture}\label{Conjecture.uniqueness}
For $n\geqslant5$ and for each admissible $i$ and $j$ where $1\leqslant i<j<n$, the
spectra~$S_{\{i,j\}_{n}^{n}}$ and~$S_{\{i,j\}_{n}^{n-1}}$ are realized by unique graphs.
\end{conjecture}
Note that Conjecture~\ref{Conjecture.uniqueness} follows from Conjecture~\ref{Conjecture.S_nnm},
the $S_{n,n}$-conjecture and Theorems~\ref{Constru.Double.Eig}, \ref{Thm.case.2.m=n-1},
and~\ref{Theorem.S.1.j.n-1.n_constr}.

\setcounter{equation}{0}
\section{Discussion}\label{section:conclusion}
Inspired by the work~\cite{FallatKirkland_et_al_2005} by Fallat et al. on the sets $S_{i,n}$, we study
the Laplacian realizability of the sets $S_{\{i,j\}_n^m}$. In the present paper, we completely described
the number of sets realizing $S_{\{i,j\}_n^m}$ for $m=n$ and $m=n-1$ and gave an algorithm for constructing
those graphs. Additionally, we also discussed some aspects of the Conjecture~\ref{Conjecture.S_nnm} about
the realizability of the set $S_{\{i,n\}_n^m}$. This conjecture correlate with the conjecture on the
realizability of sets $S_{n,n}$ introduced in~\cite{FallatKirkland_et_al_2005}.

However, there are many open problems around the Laplacian realizability of sets~$S_{\{i,j\}_n^m}$.
In the present work we fixed the number $m$ and sought for admissible numbers $i$ and $j$ absent in
Laplacian spectra of graphs. At the same time, there exist realizable sets with fixed $i$, $j$, and~$n$
with different $m$. For example, from Tables~\ref{Table.2}--\ref{Table.4} in Appendix~\ref{Tables.constructions}
it follows that the sets~$S_{\{1,3\}_{4}^{2}}$ and~$S_{\{1,3\}_{4}^{4}}$, $S_{\{2,4\}_{5}^{1}}$
and~$S_{\{2,4\}_{5}^{3}}$, $S_{\{2,4\}_{6}^{3}}$ and $S_{\{2,4\}_{6}^{5}}$ are Laplacian realizable.
Moreover, despite of Conjecture~\ref{Conjecture.uniqueness} for $m=n$ and $m=n-1$, there exist numbers $m$
such that some sets are realizable by more than one graph. For instance, the set~$S_{\{1,3\}_7^4}$
is realized by two non-isomorphic graphs $2K_{1}\vee(K_{1}\cup C_4)$ and $K_{1}\vee\overline{P_{2}\times P_{3}}$,
see Table~\ref{Table.5}. Thus, for the future work we pose the following problem.
\begin{problem}\label{problem}
Given a set $S_{\{i,j\}_{n}^{m}}$ for fixed $n$ and $i<j$, find all admissible $m$
such that $S_{\{i,j\}_{n}^{m}}$ is Laplacian realizable. Find all $m$ admitting
realizability of $S_{\{i,j\}_{n}^{m}}$ by a few non-isomorphic graphs.
\end{problem}
As well, we are interested in Conjectures~\ref{Conjecture.S_nnm}--\ref{Conjecture.uniqueness} and
will devote to them and Problem~\ref{problem} the next parts of this work.

\appendix

\setcounter{equation}{0}

\section{List of Laplacian integral graphs realizing $S_{\{i,j\}_n^m}$ up to order $7$}\label{Tables.constructions}

From~\cite{CvetkovicPetric_1984} and~\cite[p.~286]{CvetkovicRowlinson_2010} it follows that there are totally $6$ connected graphs of
order $4$, $21$ connected graphs of order $5$, and $112$ ones of order $6$. In~\cite[p.~301--304]{CvetkovicRowlinson_2010}
the authors depicted graphs up to order $5$ and found their Laplacian spectra. Thus, it is known that there are exactly
$5$ Laplacian integral graphs of order $4$ and $13$ Laplacian graphs of order $5$. In~\cite{CvetkovicPetric_1984}, the authors
depicted all these graphs of order $6$ without calculating their Laplacian spectra. In~\cite{GroneMerris_1994}, it was mentioned
that there are exactly $37$ Laplacian integral graphs. Here we list all graphs realizing sets $S_{\{i,j\}_n^m}$ up to order $6$.
We also provide some graphs of order~$7$ realizing those sets. However, finding all the graphs of order $7$ realizing $S_{\{i,j\}_n^m}$
is an open problem.

Below and throughout the text of the paper we use the following notations:
$nK_1=\bigcup\limits_{i=1}^nK_1$,
$S_n$ means the star graph on $n$ vertices, $F_n$ is the friendship graph on $2n+1$ vertices, $C_n$ is the cycle on $n$ vertices,
$P_n$ is the path graph on $n$ vertices, and $K_{p_1,\ldots,p_k}$ is the complete $k$-partite graph on $p_1+\cdots+p_k$ vertices.


%
\begin{center}
\textbf{Table 1. Laplacian integral graphs realizing $S_{\{i,j\}_n^m}$ for $n=4.$}\\
\end{center}
\begin{center}
\renewcommand*{\arraystretch}{1.9}
\begin{longtable}{|c|c|c|c}
\hline
\textbf{Construction} & \textbf{Laplacian Spectrum} & $\boldsymbol{S_{\{i,j\}_n^m}}$\label{Table.2}\\
\hline
\endfirsthead
\textbf{Construction} & \textbf{Spectrum}  \\
\hline
\endhead
\endfoot
\hline
\endlastfoot		
$S_{4}\cong K_{3,1} \cong K_{1}\vee3K_{1}$ &  ${\{0,1,1,4\}}$ & $S_{\{2,3\}_4^1}$\\
\hline
$C_{4}\cong K_{2,2}\cong 2K_{1}\vee 2K_{1}$  & ${\{0,2,2,4\}}$ & $S_{\{1,3\}_4^2}$\\
\hline
$K_{1,1,2}\cong K_{2}\vee 2K_{1}, \ \  \text{the diamond graph}$ & ${\{0,2,4,4\}}$ & $S_{\{1,3\}_4^4}$\\
\end{longtable}
\end{center}
\begin{center}
\textbf{Table 2. Laplacian integral graphs realizing $S_{\{i,j\}_n^m}$ for $n=5.$}\\
\end{center}
\begin{center}
\renewcommand*{\arraystretch}{1.9}
\begin{longtable}{|c|c|c|c}
\hline
\textbf{Construction} & \textbf{Laplacian Spectrum} & $\boldsymbol{S_{\{i,j\}_n^m}}$\label{Table.3}\\
\hline
\endfirsthead
\textbf{Construction} & \textbf{Laplacian Spectrum} & $\boldsymbol{S_{\{i,j\}_n^m}}$ \\
\hline
\endhead
\endfoot
\hline
\endlastfoot
$K_{1}\vee\overline{K_{1,1,2}}$ & $\{0,1,1,3,5\}$ & $S_{\{2,4\}_5^1}$\\
\hline
$K_{3,2}$ & $\{0,2,2,3,5\}$ & $S_{\{1,4\}_5^2}$\\
\hline
$F_2\cong K_{1}\vee 2K_{2},\ \ \text{the butterfly graph}$ &$\{0,1,3,3,5\}$ & $S_{\{2,4\}_5^3}$\\
\hline		
$K_1\vee\overline{S_4}$ & $\{0,1,4,4,5\}$ & $S_{\{2,3\}_5^4}$\\
\hline
$K_{2}\vee\overline{P_3}$ & $\{0,2,4,5,5\}$ & $S_{\{1,3\}_5^5}$\\
\end{longtable}
\end{center}

\newpage
\begin{center}
\textbf{Table 3. Laplacian integral graphs that realizing $S_{\{i,j\}_n^m}$  for $n=6$.}\\
\end{center}
\begin{center}
\renewcommand*{\arraystretch}{1.9}
\begin{longtable}{|c|c|c|c|c|c|c}
\hline
\textbf{Graph} & \textbf{Lapl. Spectrum} & $\boldsymbol{S_{\{i,j\}_n^m}}$&\textbf{Graph} & \textbf{Lapl. Spectrum} & $\boldsymbol{S_{\{i,j\}_n^m}}$\label{Table.4}\\
\hline
\endfirsthead
\textbf{Graph} & \textbf{Laplacian Spectrum}& $\boldsymbol{S_{\{i,j\}_n^m}}$  \\
\hline
\endhead
\endfoot
\hline
\endlastfoot
$K_{1}\vee(2K_{1}\cup P_3)$ &${\{0,1,1,2,4,6\}}$ & $S_{\{3,5\}_6^1}$&$K_{1}\vee\overline{K_{2,3}}$&${\{0,1,3,4,4,6\}}$ & $S_{\{2,5\}_6^4}$\\
\hline
$K_{1}\vee(K_{1}\cup S_4)$&${\{0,1,2,2,5,6\}}$ & $S_{\{3,4\}_6^2}$&$2K_{1}\vee\overline{K_{1,3}}$&${\{0,2,4,5,5,6\}}$ & $S_{\{1,3\}_6^5}$\\
\hline
$K_{1}\vee(K_{1}\cup C_4)$&${\{0,1,3,3,5,6\}}$ & $S_{\{2,4\}_6^3}$&$K_{1}\vee(K_{1}\cup K_{1,1,2})$ & ${\{0,1,3,5,5,6\}}$ & $S_{\{2,4\}_6^5}$\\
\hline
$\overline{P_{2}\times P_{3}}$&${\{0,1,3,3,4,5\}}$ & $S_{\{2,6\}_6^3}$&$P_{3}\vee\overline{P_3}$&${\{0,3,4,5,6,6\}}$ & $S_{\{1,2\}_6^6}$\\
\hline
$P_{2}\times P_{3}$&${\{0,1,2,3,3,5\}}$ & $S_{\{4,6\}_6^3}$&$K_{2}\vee\left(K_{1}\cup P_3\right)$ & ${\{0,2,3,5,6,6\}}$ & $S_{\{1,4\}_6^6}$\\
\end{longtable}
\end{center}
%

\begin{center}
\textbf{Table 4. Some Laplacian integral graphs that realizing $S_{\{i,j\}_n^m}$  for $n=7$.}\\
\end{center}
\begin{center}
\renewcommand*{\arraystretch}{1.9}
\begin{longtable}{|c|c|c|c}
\hline
\textbf{Construction} & \textbf{Laplacian Spectrum} & $\boldsymbol{S_{\{i,j\}_n^m}}$\label{Table.5}\\
\hline
\endfirsthead
\textbf{Construction} & \textbf{Laplacian Spectrum} & $\boldsymbol{S_{\{i,j\}_n^m}}$  \\
\hline
\endhead
\endfoot
\hline
\endlastfoot
$K_{1}\vee(P_3\cup \overline{P_3})$&${\{0,1,1,2,3,4,7\}}$ & $S_{\{5,6\}_7^1}$\\
\hline
$K_{1}\vee[2K_{1}\cup(K_{1}\vee\overline{P_3})]$&${\{0,1,1,2,4,5,7\}}$ & $S_{\{3,6\}_7^1}$\\
\hline	
$2K_{1}\vee(K_{1}\cup S_4)$&${\{0,2,3,3,5,6,7\}}$ & $S_{\{1,4\}_7^3}$\\
\hline
$2K_{1}\vee(K_{1}\cup C_4)$ & ${\{0,2,4,4,5,6,7\}}$ & $S_{\{1,3\}_7^4}$\\
\hline
$K_{1}\vee\overline{P_{2}\times P_{3}}$ & ${\{0,2,4,4,5,6,7\}}$ & $S_{\{1,3\}_7^4}$\\
\hline
$K_{1}\vee(P_{2}\times P_{3})$ & ${\{0,2,3,4,4,6,7\}}$ & $S_{\{1,5\}_7^4}$\\
\hline
$2K_{1}\vee(K_{2}\cup P_3)$&${\{0,2,3,4,5,5,7\}}$ & $S_{\{1,5\}_7^5}$\\
\hline
$K_{1}\vee[K_{1}\cup(K_{1}\vee\overline{S_4})]$&${\{0,1,2,5,5,6,7\}}$ & $S_{\{3,4\}_7^5}$\\
\hline
$2K_{1}\vee(K_{1}\cup K_{1,1,2})$&${\{0,2,4,5,6,6,7\}}$ & $S_{\{1,3\}_7^6}$\\
\hline
$K_{1}\vee [K_{1}\cup (K_{2}\vee\overline{P_3})]$&${\{0,1,3,5,6,6,7\}}$ & $S_{\{2,4\}_7^6}$\\
\hline
$K_{1}\vee[2K_1\vee (K_{1}\cup P_3)]$&${\{0,3,4,5,6,7,7\}}$ & $S_{\{1,2\}_7^7}$\\
\hline
$K_{2}\vee[K_1\cup(K_1\vee\overline{P_3})]$&${\{0,2,3,5,6,7,7\}}$ & $S_{\{1,4\}_7^7}$\\
\hline
$K_{2}\vee (K_{2}\cup P_3)$&${\{0,2,3,4,5,7,7\}}$ & $S_{\{1,6\}_7^7}$\\
\end{longtable}
\end{center}	
%

\section*{Acknowledgement}
The work of M.\,Tyaglov was partially supported by National Natural Science Foundation of China under grant no.~11871336.

\bibliographystyle{abbrv}
\bibliography{main.bib}

\end{document}